
\magnification=\magstep1
\def\to{\ \longrightarrow\ }

\def\nl{\hfill\break}

\def\hexnumber#1{\ifcase#1 0\or 1\or 2\or 3\or 4\or 5\or 6\or 7\or 8\or
 9\or A\or B\or C\or D\or E\or F\fi}
%
%
\font\twelvemsa=msam10 scaled 1200   
\font\tenmsa=msam10                  
\font\ninemsa=msam9            \font\sevenmsa=msam7
\font\sixmsa=msam6             \font\fivemsa=msam5
%
%
\newfam\msafam                 \textfont\msafam=\tenmsa
\scriptfont\msafam=\sevenmsa   \scriptscriptfont\msafam=\fivemsa
\edef\hexa{\hexnumber\msafam}        
\def\msa{\fam\msafam\tenmsa}         
%
%
\font\twelvemsb=msbm10 scaled 1200   
\font\tenmsb=msbm10                  
\font\ninemsb=msbm9            \font\sevenmsb=msbm7
\font\sixmsb=msbm6             \font\fivemsb=msbm5
%
\newfam\msbfam                 \textfont\msbfam=\tenmsb       
\scriptfont\msbfam=\sevenmsb   \scriptscriptfont\msbfam=\fivemsb
\edef\hexb{\hexnumber\msbfam}        
\def\msb{\fam\msbfam\tenmsb}         
%
%
\font\twelveeufm=eufm10 scaled 1200  
\font\teneufm=eufm10                 
\font\nineeufm=eufm9           \font\seveneufm=eufm7
\font\sixeufm=eufm6            \font\fiveeufm=eufm5
%
\newfam\eufmfam                \textfont\eufmfam=\teneufm
\scriptfont\eufmfam=\seveneufm \scriptscriptfont\eufmfam=\fiveeufm
\edef\hexf{\hexnumber\eufmfam}      
\def\frak{\fam\eufmfam\teneufm}     
%
%
%
\font\twelverm=cmr10 scaled 1200    
\font\ninerm=cmr9                   
\font\sixrm=cmr6   
%
\font\twelvei=cmmi10 scaled 1200    
\font\ninei=cmmi9                   
\font\sixi=cmmi6  
%
\font\twelvesy=cmsy10 scaled 1200   
\font\ninesy=cmsy9                  
\font\sixsy=cmsy6  
%
\font\twelvebf=cmbx10 scaled 1200   
\font\ninebf=cmbx9                  
\font\sixbf=cmbx6  
%
%
\font\twelveit=cmti10 scaled 1200   
\font\nineit=cmti9                  
%
\font\twelvesl=cmsl10 scaled 1200   
\font\ninesl=cmsl9                  
%
\font\twelvett=cmtt10 scaled 1200   
\font\ninett=cmtt9                  
%
%
%
%
\def\small{%
%
%
\textfont0=\ninerm \scriptfont0=\sixrm \scriptscriptfont0=\fiverm
\def\rm{\fam0\ninerm}        
%
%
\textfont1=\ninei \scriptfont1=\sixi \scriptscriptfont1=\fivei
%
%
\textfont2=\ninesy \scriptfont2=\sixsy \scriptscriptfont2=\fivesy
%
%
\textfont3=\tenex \scriptfont3=\tenex \scriptscriptfont3=\tenex
%
%
\textfont\bffam=\ninebf \scriptfont\bffam=\sixbf
\scriptscriptfont\bffam=\fivebf \def\bf{\fam\bffam\ninebf}%
%
%
\textfont\itfam=\nineit \def\it{\fam\itfam\nineit}%
\textfont\slfam=\ninesl \def\sl{\fam\slfam\ninesl}%
\textfont\ttfam=\ninett \def\tt{\fam\ttfam\ninett}%
%
%
%
\textfont\msafam=\ninemsa \scriptfont\msafam=\sixmsa
\scriptscriptfont\msafam=\fivemsa \def\msa{\fam\msafam\ninemsa}%
%
%
\textfont\msbfam=\ninemsb \scriptfont\msbfam=\sixmsb
\scriptscriptfont\msbfam=\fivemsb \def\msb{\fam\msbfam\ninemsb}%
%
%
\textfont\eufmfam=\nineeufm  \scriptfont\eufmfam=\sixeufm
\scriptscriptfont\eufmfam=\fiveeufm \def\frak{\fam\eufmfam\nineeufm}%
%
%
%
\normalbaselineskip=11pt
\setbox\strutbox=\hbox{\vrule height8pt depth3pt width0pt}%
%
%
\normalbaselines\rm}    
%
%
%
%
\def\large{%
\textfont0=\twelverm \scriptfont0=\ninerm \scriptscriptfont0=\sevenrm
\def\rm{\fam0\twelverm}%
\textfont1=\twelvei \scriptfont1=\ninei \scriptscriptfont1=\seveni
\textfont2=\twelvesy \scriptfont2=\ninesy \scriptscriptfont2=\sevensy
\textfont3=\tenex \scriptfont3=\tenex \scriptscriptfont3=\tenex
\textfont\bffam=\twelvebf \scriptfont\bffam=\ninebf
\scriptscriptfont\bffam=\sevenbf \def\bf{\fam\bffam\twelvebf}%
\textfont\itfam=\twelveit \def\it{\fam\itfam\twelveit}%
\textfont\slfam=\twelvesl \def\sl{\fam\slfam\twelvesl}%
\textfont\ttfam=\twelvett \def\tt{\fam\ttfam\twelvett}%
\textfont\msafam=\twelvemsa \scriptfont\msafam=\ninemsa
\scriptscriptfont\msafam=\sevenmsa \def\msa{\fam\msafam\twelvemsa}         
\textfont\msbfam=\twelvemsb \scriptfont\msbfam=\ninemsb
\scriptscriptfont\msbfam=\sevenmsb \def\msb{\fam\msbfam\twelvemsb}         
\textfont\eufmfam=\twelveeufm  \scriptfont\eufmfam=\nineeufm
\scriptscriptfont\eufmfam=\seveneufm \def\frak{\fam\eufmfam\teneufm}
\normalbaselineskip=15pt
\setbox\strutbox=\hbox{\vrule height11pt depth4pt width0pt}%
\normalbaselines\rm}%
%
\def\Bbb{\msb}

%

%
\mathchardef\plussquare="0\hexa01
\mathchardef\nge="3\hexb0B
\mathchardef\maltesecross="0\hexa7A
\mathchardef\del="0\hexf01
%

%

\input epsf
\overfullrule=0pt

\font\npt=cmr9
\font\Bbb=msbm10

\font\secfont=cmbx10
\font\ab=cmbx8

\font\nam=cmr8
\font\aff=cmti8

\font\em=cmti10

\mathchardef\square="0\hexa03
\def\qed{\hfill$\square$\par\rm}
\def\np{\vfill\eject}
\def\boxing#1{\ \lower 3.5pt\vbox{\vskip 3.5pt\hrule \hbox{\strut\vrule
\ #1 \vrule} \hrule} }

\def\down#1{\ \lower 3.5pt\vbox{\vskip 3.5pt \hbox{\strut \ #1 \vrule} \hrule} }
\def\negdown#1{\ \lower 3.5pt\vbox{\vskip 3.5pt \hbox{\strut  \vrule \ #1 }\hrule} }

\hsize=6.3 truein
\vsize=9 truein

\baselineskip=13 pt
\parskip=\baselineskip
 1

\parindent=0pt

\def\Z{\hbox{\Bbb Z}}
\def\I{\hbox{\Bbb I}}
\def\R{\hbox{\Bbb R}}
\def\C{\hbox{\Bbb C}}
\def\H{\hbox{\Bbb H}}

\def\op{\buildrel o \over p}
\def\np{\buildrel n \over p}
\def\oM{\buildrel o \over {\cal M}}
\def\oP{\buildrel o \over {P}}

\def\nM{\buildrel n \over {\cal M}}
\def\nP{\buildrel n \over {P}}

\def\oover#1{\vbox{\ialign{##\crcr
{\npt o}\crcr\noalign{\kern 1pt\nointerlineskip}
$\hfil\displaystyle{#1}\hfil$\crcr}}}







\newif \iftitlepage \titlepagetrue

\def\diagram{\global\advance\diagramnumber by 1
$$\epsfbox{treffig.\number\diagramnumber}$$}
\def\ddiagram{\global\advance\diagramnumber by 1
\epsfbox{treffig.\number\diagramnumber}}

\newcount\diagramnumber
\diagramnumber=0

\newcount\secnum \secnum=0
\newcount\subsecnum
\newcount\defnum
\def\section#1{
                \vskip 10 pt
                \advance\secnum by 1 \subsecnum=0
                \leftline{\secfont \the\secnum \bf\quad#1}
                }

\def\subsection#1{
                \vskip 10 pt
                \advance\subsecnum by 1 
                \defnum=1
                \leftline{\secfont \the\secnum.\the\subsecnum\ \rm\quad #1}
                }

\def\definition{
                \advance\defnum by 1 
                \bf Definition 
\the\secnum .\the\defnum \rm \ 
                }

\def\lemma{
                \advance\defnum by 1 
                \par\bf Lemma  \the\secnum
.\the\defnum \rm \ \par
                }

\def\theorem#1{
                \advance\defnum by 1 
                \par\bf Theorem  \the\secnum
.\the\defnum \sl \ #1 \par\rm
               }

\def\cite#1{
				\secfont [#1]
				\rm$\!\!\!$\nobreak
}

\vglue 20 pt
\centerline{\secfont Tackling the Trefoils}
\medskip

\centerline{\nam Roger Fenn\footnote{${}^1$}{\aff School of Mathematical Sciences, University of Sussex, Falmer, Brighton, BN1 9RH, England, e-mail address: rogerf@sussex.ac.uk}}

\bigskip
\centerline{\nam ABSTRACT}
\leftskip=0.25 in
\rightskip=0.25in


{\ab The classical trefoil is famous for having a three-colouring which distinguishes it from the unknot. The three-colouring is also notorious for not distinguishing the right handed from the left handed trefoil. However with a bit of tweaking the three colours can also be used for this task. What lies behind the method is a new operation on biracks called doubling which converts the 3-colour quandle into a biquandle. Colouring with this biquandle distinguishes the right handed from the left handed trefoil. Equivalently it defines an element of the homology of the quandle or biquandle classifying space.}

\leftskip=0 in
\rightskip=0in
{\npt Keywords: 3-colouring, virtual knots, trefoil, quandles, racks, biquandles and their homology.}
\section{\bf Introduction}
There have been many proofs that the right and left handed trefoils are distinct. 
 The earliest proof that the trefoil is inequivalent to its mirror image is due to Max Dehn, who proved it by tracking the longitude in the fundamental group. There is a good exposition of this and reference to Max Dehn in the book \cite{S}.
Perhaps the simplest proof  is via the  Kauffman bracket polynomial, starting the construction from scratch and using the Reidemeister moves. The next simplest to that is using the third Vassiliev invariant. 

What is interesting about these and other  solutions  to different problems is how they drive forward the search and invention of new mathematical methods and engines. In this case it is the algebraic topology related to the classifying space of a quandle or biquandle. Indeed quandles and racks were invented in the search for classical knot invariants whilst biquandles turned out to be very useful for virtual knot theory.

Pictured is a classical trefoil (both right and left handed versions), with a 3-colouring using the colours red (r), green (g), and blue (b) which shows that they are non-trivial, but fails to distinguish right from left.
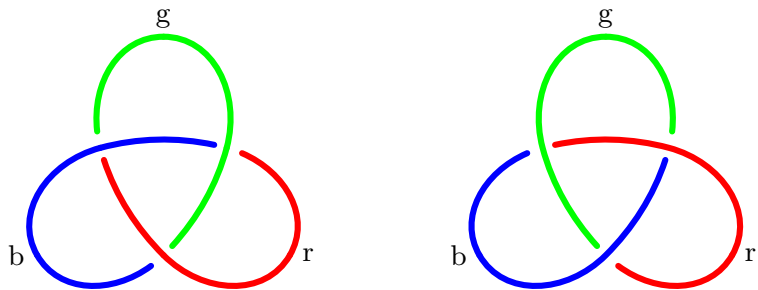
\diagram
However the next picture shows a different 3-colouring using pairs of colours; red/green (rg), blue/red (br), green/green, (gg), red/red (rr) and red/blue (rb) which does distinguish the right trefoil from the left. 
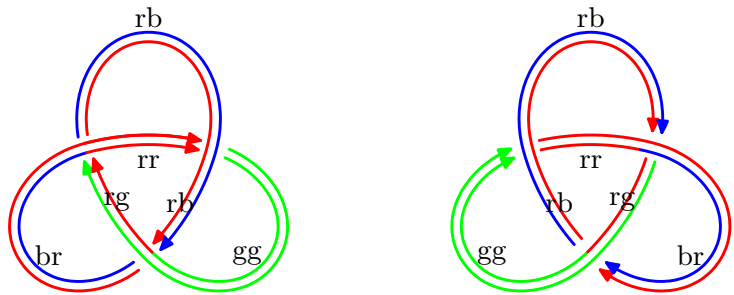
\diagram
In this paper it will be shown why this is the case and how it relates to homology of quandles and biquandles.

I would like to thank Lou Kauffman, Jozef Przytycki, Vassily Manturov, Colin Rourke and Brian Sanderson for helpful comments. I would particularly like to thank Seiichi Kamada for pointing out amongst other things that the first draft contained an ambiguity with orientation conventions.
\section{\bf Colouring by pairs}
The edges of a diagram are coloured with pairs of the three colours, red, green and blue.
The pairs of colours can change on the overcrossing as well as the undercrossing. Moreover the order is important. For example (gr) means that red is on the left and green is on the right as one traverses the diagram in the given direction.
\diagram
Not all pairs of colours are allowed at a crossing. The possibilities are indicated in the following diagram.
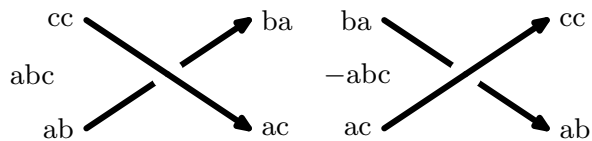
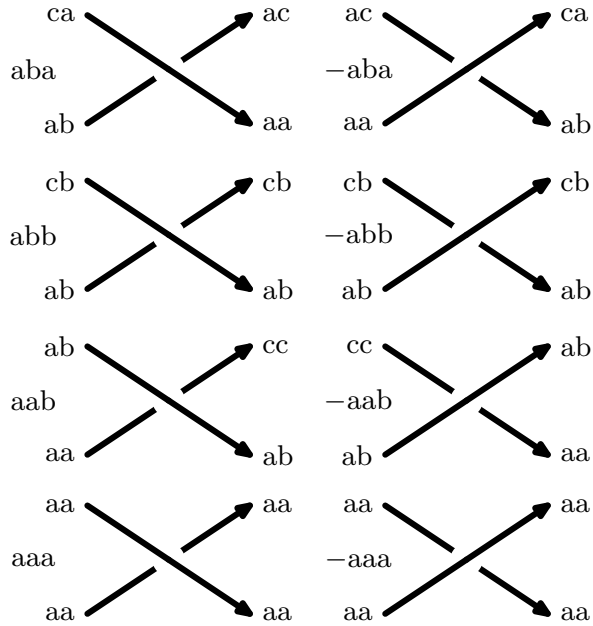
\diagram\diagram\diagram\diagram\diagram
The letters $abc$ are a permutation of the colours $rgb$. There are 27 possible colourings of the positive crossings. They are labelled $abc$, $aba$, $abb$, $aab$ and $aaa$, for reasons which will become clear later. 

The 27 possible colourings of negative crossings are shown on the right. They can be obtained by reflecting the positive crossings in a vertical line.

If the diagram contains a virtual crossing then the pairs cross unaltered as follows.
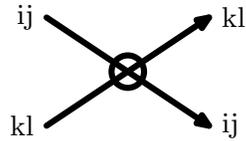
\diagram
The letters $ijkl$ are combinations of the colours $rgb$ possibly with repitition. 
\section{\bf Reidemeister moves}

We now see how the colouring by pairs is affected during a Reidemeister move.

Suppose an edge is coloured $ab$ and a curl (monogon) is introduced by an expansive Reidemeister move of type I. Then a crossing of type $\pm abb$  is introduced. If the edge is coloured $aa$ then a crossing of type $\pm aaa$  is introduced. Conversely only crossings of type $\pm abb$ or $\pm aaa$ can be the vertices in a curl and the curl can be removed if desired. Recall that $abc$ is a permutation of the original three colours.

A similar description applies to a virtual Reidemeister move of type I. 

A positive crossing can be cancelled with a negative crossing by a Reidemeister move of type II since any pair of colours is matched by a unique pair of colours with an opposite crossing either on the sides or top/bottom.

On the other hand an expansive Reidemeister move of type II presents difficulties since the form of the pairs have to match the limited number of colourings of the crossings. However without virtual crossings this does not present a problem as we shall now see. Even with the presence of virtual crossings there are situations where expansive Reidemeister moves of type II can take place.

Let $D$ be a classical knot diagram (or irreducible link) and orient all the regions anti-clockwise. An edge, $e$, on the boundary of a region $F$ may or may not be oriented coherently with $F$. Suppose the edge is oriented coherently and is coloured by the pair $ab$, where $a\ne b$ then colour the region $F$ by $c$. If the edge is not oriented coherently and is coloured by the pair $ab$, colour the region by $a$. If $e$ is coloured by $aa$ then colour $F$ by $a$.

\theorem{With the notation above, the colour of the region $F$ is independent of the boundary edge $e$ chosen.}
{\bf Proof} Let $e_1$ and $e_2$ be two edges on the boundary of $F$. By induction we may as well assume that they have a crossing in common. The result follows after a consideration of all the possibilities. \qed

There are 45 colouring possibilities for the Reidemeister III move and we shall not consider all of them since we shall see the general pattern later in the paper. However to get an idea look at the following pair coloured Reidemeister III move.
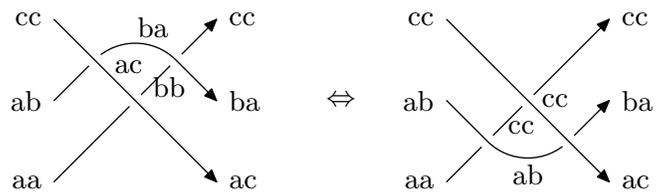
\diagram
\section{\bf Crossing Invariants}
To get an invariant we sum up the crossing points of the diagram with sign. So for example the right hand trefoil gives $+rbg+rrb+rgr$. The left hand trefoil gives the negative,  $-rbg-rrb-rgr$

In order to obtain an invariant we look at how this sum changes under the Reidemeister moves. During a Reidemeister move of type I, crossing points of the form $abb$ and $aaa$ can be eliminated or constructed at will so these will be put equal to zero. For reasons which will be clear later also put $aab=0$.

Reidemeister moves of type II confirm the inverse of a crossing.

Reidemeister moves of type III imply relations amongst the crossings. For example, the Reidemeister move of type III illustrated above implies the relation $abc+aac+bba=aab+ccc+abc$, confirming the r\^oles of $aac, bba, aab$ and $ccc$ as zeros.

It turns out that the element $abc+aca$ is of order 3 in the resulting abelian group. This shows that the two trefoils are distinct.
\section{Two sided knot diagrams}
In an oriented classical knot diagram each edge has a left and right adjacent region. Moreover these two regions are distinct. This may not happen with a virtual knot diagram. If we think of a virtual knot as a knot on a surface, then a meridian on a torus divides the surface into one region. Alternatively consider the {\bf virtual trefoil} in the following figure.
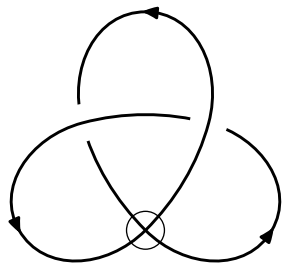
\diagram

We now look at conditions which allow a diagram to be two sided. Since over and undercrossings are irrelevant to this problem we may consider only {\bf flat virtual knots}. They are represented by an oriented 4-valent graph in a surface which is the image in general position of an immersed circle or a number of circles. Each crossing is a vertex and the edges inherit an orientation, so that at each vertex they enter and leave on opposite sides. Two such diagrams are equivalent if they are related by a sequence of homeomorphisms, local Reidemeister type moves and surgeries disjoint from the image. A diagram is {\bf minimal} if the surface in which it lies has minimal genus. 
The {\bf genus} of a virtual knot is this minimal genus. So classical knots have genus zero. 

By a result of Kuperberg, \cite{Kuper}, two minimal diagrams representing the same class are related by a sequence of homeomorphisms and local Reidemeister type moves. Note that if a diagram is minimal and connected then the regions are open discs.

A virtual knot, $K$, flat or otherwise, is called {\bf 2-sided} if it has a diagram so that any edge is on the boundary of exactly 2 regions. It is said to be {\bf irreducible} if at every crossing there are exactly 4 regions. A minimal irreducible 2-sided diagram is called {\bf cellular}.

A cellular diagram is the 1-skeleton of a regular cellulation of the surface. We say that a cellular diagram is {\bf 2-colourable} if the regions can be coloured, chess board fashion, with 2 colours. This is a generalisation of a more general colouring with a biquandle which is defined later in the paper. For example on a torus take a knot with 2 components, a meridian and a longitude. If we double each component then we have a cellular diagram whose regions can be 2-coloured. If we triple the meridian and longitude we have a cellular diagram which cannot be 2-coloured.

Cellular diagrams which are 2-colourable are called {\bf atoms}, see \cite{B}.

We can decide which cellular diagram can be 2-coloured by looking at the orientation of the edges. Colour the edges alternately black and white along the diagram. This is always possible if the diagram represents a knot with one component but may not be possible for knots with more than one component. If an edge is coloured white keep its natural orientation. If it is black, reverse its orientation. Then each crossing is a source, sink or saddle.
\diagram
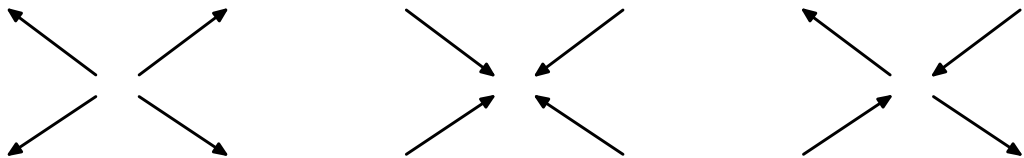
\centerline{\npt Source, sink and saddle}\smallskip
The property of being alternately oriented is easily seen to be invariant under the Reidemeister moves and so is a knot invariant.

If the knot has one component then there is a convenient way of seeing the crossing type by means of the {chord diagram}. A crossing is an {\bf odd crossing} if its chord crosses an odd number of other chords. Otherwise it is {\bf even}. Odd crossings correspond to sources or sinks; even crossings to saddles.

If all crossings are saddles then the  alternately oriented diagram is called {\bf good}.  This is similar to the treatment in \cite{IKK, YM} where it is called a magnetic graph.  
\theorem{A cellular knot diagram has a 2-colouring if and only if it is a good alternately oriented diagram}
{\bf Proof}
If a cellular knot diagram has a 2-colouring and if the underlying surface is oriented, then orient the boundary of each black cell according to the orientation of the surface. This defines a good alternately oriented diagram. Conversely a good alternately oriented minimal diagram defines a 2-colouring of the cells by the reverse process.
\qed
\diagram
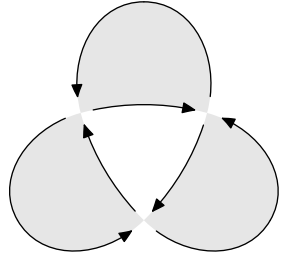
\centerline{\npt Good alternatively oriented trefoil with 2-colouring}\smallskip
\diagram
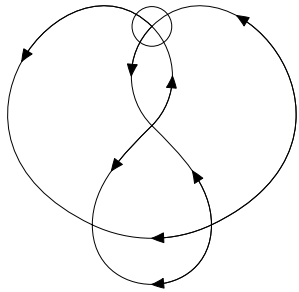
\centerline{\npt Alternatively oriented virtual diagram with sink, source and saddle}\smallskip
\theorem{Let $D$ be an alternatively oriented virtual diagram. Then
$$\#\hbox{sinks}=\#\hbox{sources}.$$}
{\bf Proof}  Consider the chord diagram of $D$. This consists of $n$ circles corresponding to the $n$ components and two types of chords. The interior chords join points on the same circle and correspond to self crossings. The exterior chords join different circles and correspond to where they cross. Note that a necessary and sufficient condition for the diagram to be alternately oriented is that there are an even number of exterior chords attached to each circle.

We can simplify the diagram without changing the conclusion of the result by interchanging end points of adjacent internal and external chords. In this manner all the external chord end points can be grouped together as can those of the internal chords.

Consider the internal chords. Define cycles as follows. The vertices are the chords and the edges are the pieces of the circles at the end points travelling anticlockwise to the adjacent chord end. Around each cycle the edges are oriented by arrows in one of two ways. It is at a change of orientation that a sink or source is generated. If the arrows converge at the vertex then it is a sink. Otherwise it is a source. A simple counting argument shows that there is an equal number of sinks and sources.

A similar argument now works for the exterior chords. \qed

Let $K$ be a knot which has an alternately oriented diagram and let  $\chi(K)$ be the minimum number of sinks (sources) for any diagram representing $K$. Note that $\chi=0$  is a necessary and sufficient condition for $K$ to be represented by a 2-coloured diagram.

The virtual trefoil shown earlier has one source and one sink. So in this case $\chi=1$. 

I am indebted to Andrew Bartholomew for the following example of a virtual knot with a good alternate orientation.
\diagram
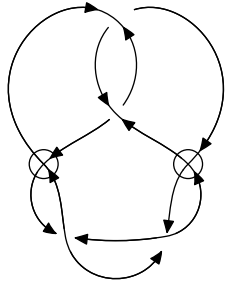
\centerline{\npt Good alternatively oriented virtual knot}\smallskip
The usual representation of virtual knot diagrams means that this corresponds to a diagram on an orientable surface of genus 2. However if we thicken the diagram then it lies on a thickened 4-valent graph with 4 vertices and 8 edges which is the image of an immersed ribbon. The diagram is 2-sided and the surface has genus 1 and 4 boundary components. Of course, to show that this is the minimal genus we would have to show that it is not classical. 
\diagram
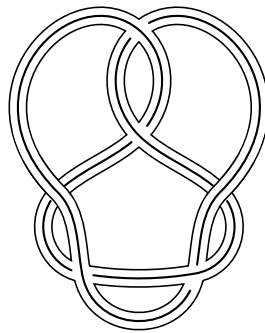
\centerline{\npt Virtual knot on a torus with 4 holes}\smallskip
\section{\bf Racks, biquandles etc: }
We will slightly extend the definition of a birack and a biquandle to suit our needs later in the paper. The more usual notation and definitions can be found in \cite{FJK}, or \cite{BaF} where a list of small biracks etc can be found.

Let $X$ be a set; the colouring or labelling set. Let $Y$ be a subset of $X^2$.  A function $S:Y\to Y$ defines two binary operations by the formula 
$$S(a,b)=(b^a,a_b)$$
The two binary operations,
$$(a,b)\to a^b\hbox{ and }(a,b)\to a_b$$
are called  {\bf up} and {\bf down}, respectively. In previous treatments, $Y=X^2$, so the binary operations are defined for all $a,b\in X$. This extension of the definition means that $b^a$ and $a_b$ are only defined if $(a,b)\in Y$ and will give us the flexibility needed later. However whenever we write  $b^a$ and $a_b$ we will always assume that the operations are defined.

The convenience of the exponential and suffix notation is that brackets can be inserted in an obvious fashion and so are not needed, For example
$$a^{bc}=(a^b)^c,\ a^{b_c}=a^{(b_c)},\ {a^b}_c=(a^b)_c, \hbox{ etc. }$$
On the other hand expressions such as $a^b_c$ are ambiguous and are not used.

Think of $b$ as acting on $a$ in both cases. We want these actions to be invertible. So there are binary operations
$$(a,b)\to a^{b^{-1}}\hbox{ and }(a,b)\to a_{b^{-1}}$$
satisfying
$$a^{bb^{-1}}=a^{b^{-1}b}=a\hbox{ and }a_{bb^{-1}}=a_{b^{-1}b}=a$$
So each suitable element $a$ of $X$ defines two permutations given by $x\to x^a$ and $x\to x_a$.

It is convenient at this stage to introduce the function $F$ defined by $F(a,b^a)=(b, a_b)$.  In \cite{FJK} this was called the sideways map, $S_{-}^{+}$.

After these preliminaries we list the three axioms, B1-3 needed to define a {\bf biquandle}.

{\bf B1: } $F$ is invertible and preserves the diagonal, $\{(a,a)|a\in X\}$.

{\bf B2: } $S$ is invertible. So there is a function $\overline{S}: X^2\to X^2$ such that $S\overline{S}=\overline{S}S=id$

{\bf B3: } {Let $S_1=S\times id$ and $S_2=id\times S$ then $S_1S_2S_1=S_2S_1S_2$}

These axioms are a consequence of the three Reidemeister moves. If only B2 and B3 are satisfied then we have a {\bf birack}.

The function $F$ satisfies
$$F(a,a)=(a^{a^{-1}},a_{a^{a^{-1}}}).$$
Axiom B1 implies that for all $a \in X$ there is a unique $x\in X$ such that $a^x=x, a=x_a$  and there is a unique $y\in X$ such that $a_y=y, a=y^a$.
We have $x=a_{a^{-1}}$, $y=a^{a^{-1}}$ and so B1 is equivalent to
$a^{a_{a^{-1}}}=a_{a^{-1}}\hbox{ and }a_{a^{a^{-1}}}=a^{a^{-1}}$.

In \cite{Stan} it is shown that only one half of B1 is necessary.

An amusing implication of B1 is that the infinite tower and the infinite well
$$x=a^{a^{a^{.^{.^{.^.}}}}},\quad y=a_{a_{a_{._{._{._.}}}}}$$
make sense.

Using B2 write $\overline{S}(a,b)=(b_{\overline{a}},a^{\overline{b}})$.
This defines two binary operations
$$(a,b)\to a^{\overline b}\hbox{ and }(a,b)\to a_{\overline b}$$
In terms of the previous operations $\overline{a^b}$ acts down on $b$ as $a^{-1}$ and $\overline{a_b}$ acts up on $b$ as $a^{-1}$. So all these new operations are defined by the two initial up and down operations.

B3 is sometimes called the set theoretic Yang-Baxter equation. If we follow the progress of the triple $(a,b,c)$ through the two sides of the equation and swap variables we arrive at three relations true for all $a,b,c\in X$.
$$a^{c_bb^c}=a^{bc},\quad {a^b}_{c^{b_a}}={a_c}^{b_{c^a}},\quad a_{c^bb_c}=a_{bc}$$

If the down operation is trivial, so $a_b=a$ for all $a,b\in X$, then a biquandle becomes a quandle. Symmetrically, this is also the case if the up operation is trivial. A birack with trivial down (up) operation is a rack. 
\section{\bf Examples}
For many examples of racks and quandles, see \cite{FR}.

The simplest quandle, with both up and down operations trivial, is the {\bf twist}. In this case $S(a,b)=(b,a)$. If $X$ has $n$ elements denote the twist by $\I_n$. 

The {\bf 3-colour} quandle is the reflection set of transpositions in the symmetric group $S_3$. The quandle operation is conjugation. It is notated $Q^3_3$ in \cite{BaF}. It is related to the three colouring given in the introduction. Its double is the colouring by pairs.

The {\bf black and white} biquandle is the smallest biquandle which is not a rack. It has 2 elements $\{b, w\}$ and each permutation is the transposition $(bw)$. It is notated $BQ^2_1$ in \cite{BaF} and is an important example related to 2-colouring.

If a biquandle  has the property that the up operation (or the down operation) on its own defines a quandle then it is called a {\bf quandle-related} biquandle associated to that quandle.  In a similar fashion we can define a {\bf rack related} birack.

Our generalisation of biracks given earlier allows us to define the {\bf double} of a birack. Let $X$ be a birack and let $Z$ be the set of pairs of pairs $\{(a,b), (c,d)| a, b, c, d\in X\}$. Then $Y$ is the subset 
$$Y=\{(a^b,c),(a,b) | a,b,c,\in X\}.$$ 
The doubled operations are 
$$(a,b)^{(a^b,c)}=(a^{c_b},b^c)\ \hbox{ and }\ (a^b,c)_{(a,b)}=(a,c_b)$$
Doubling converts racks into biracks and quandles into biquandles.

For example consider the double of the 3-colour quandle. This is related to the colouring of pairs on the right and left trefoil.

A biquandle is said to be {\bf linear} if it is determined by a $2\times2$ matrix $S=\pmatrix{A&B\cr C&D\cr}$, where $A, B, C, D$ are elements of an associative ring related by certain equations, see for example \cite{BuF, BuF2, FT}. The elements $A, B, C, D$ satisfy

$${\cal F}:\ A^{-1}B^{-1}AB-B^{-1}AB=BA^{-1}B^{-1}A-A$$
and $C, D$ are defined by
$$C=A^{-1}B^{-1}A(1-A),\quad D=1-A^{-1}B^{-1}AB.$$

The only example of a commutative linear biquandle is given by
$$a^b=\lambda a+(1-\lambda\mu)b,\quad a_b=\mu a$$
where $a,b\in X$, a $\Z[\lambda^{\pm1}, \mu^{\pm1}]$-module \cite{Swa}. This is called the {\bf Alexander} biquandle and is denoted by $A_{\lambda\mu}(X)$. If $\mu=1$ then resulting quandle is called the {\bf Burau} quandle and is denoted by $B_\lambda(X)$. If $\mu=\lambda=1$ then we have the twist.

There are many examples of linear non-commutative biquandles.  For example let $\H$ denote the quaternion algebra with standard generators $1,i, j, k$. If $X$ is a left $\H$-module then 
$$a^b=ja+(1+i)b,\quad a_b=-ja+(1+i)b$$
is a biquandle. This is called the {\bf Budapest} biquandle and is just one of a huge family of linear non-commutative biquandles described in  \cite{F, BF, BuF, BuF2, FT}.



\section{Colouring a virtual knot by biracks}
Let $D$ be a diagram either on a surface $\Sigma$ or in the plane with virtual crossings and let $X$ be a birack. An {\bf edge colouring} of $D$ by $X$ is an assignment of each edge to an element of $X$, its {\bf colour}, such that at each crossing, positive, negative and virtual, the colours, $b, c, b^c, c_b\in X$ satisfy the conditions illustrated in the following figure.
\diagram
\centerline{\npt Edge colouring}\smallskip
For example the first figure of this paper shows the trefoil coloured by the 3-colouring quandle. The second figure shows the trefoil coloured by the double of the 3-colouring.
An alternatively oriented diagram has been edge coloured by the black and white biquandle.

In a {\bf whole-colouring} of $D$ by $X$ the regions of the diagram are also labelled so that at any edge the left and right are labelled as follows. 
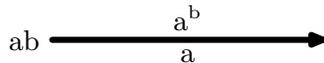
\diagram
The edge is now labelled by the pair $ab$, where (perversely) the region on the right is labelled $a$ and the region on the left is labelled $a^b$. The residual labelling of the edge is by $b$. Of course if $a\ne a^b$ then the diagram has to be 2-sided.

Note that any colouring of a region transmits around the diagram and determines the colouring of all the other regions.
there is also a dual colouring in which the action of the edge colouring on neighbouring regions is by the down operation.

At crossings the whole-colouring looks as follows.
\diagram
\centerline{\npt Whole colouring}

In the example at the beginning of the two kinds of trefoil, the colouring is now by the double of the 3-colour quandle.

An edge of a 2-sided diagram coloured by a birack can be converted into an edge coloured by the double.

To determine whether a virtual knot diagram is 2-sided let us suppose that the edges are coloured alternatively by $\{b, w\}$, that is by  $BQ^2_1$. If this is possible it can be done in two ways if the knot has one component. At any crossing the colours are cyclically ordered $b, b, w, w$. Now try and colour the regions. At a saddle the edges are now coloured by the pairs $bb, bb, ww, ww$ or $bw, bw, wb, wb$ in cyclic order.  However at a sink say, the pairs are $ww, wb, bb, bw$ in cyclic order.

\section{The homology of racks and biracks}
Unusually in the history of mathematics, the discovery of the homology and classifying space of a rack can be precisely dated to 2 April 1990, see \cite{history}.

Let $X$ be a birack. Associated with each $n$-tuple $(x_1, \dots, x_n)$ of elements of $X$ is a cube of dimension $n$. Each cube is canonically identified with the standard $n$-cube and they fit together to make a classifying space, see \cite{FRSd}.

The homology of this space can be calculated as follows.
 Let $C_n^{\rm BR}(X)$ be the free  abelian group generated by the
$n$-tuples $(x_1, \dots, x_n)$. Define a homomorphism
$\partial_{n}: C_{n}^{\rm BR}(X) \to C_{n-1}^{\rm BR}(X)$ by 
$$
\matrix{%
\partial_{n}(x_1, x_2, \dots, x_n)  =&
\sum_{i=2}^{n} (-1)^{i}\left[ (x_1, x_2, \dots, x_{i-1}, x_{i+1},\dots, x_n)\right. \cr
&\qquad- (x_1 ^{x_i}, x_2^{x_i}, \dots, x_{i-1}^{x_i}, {x_{i+1}}_{x_i}, \dots, {x_{n}}_{x_i}) ]\cr}
$$
for $n \geq 2$  and $\partial_n=0$ for $n \leq 1$. 
 Then
$C_\ast^{\rm BR}(X)
= \{C_n^{\rm BR}(X), \partial_n \}$ is a chain complex and consequently has homology groups with any coefficient groups.

For example consider the black and white biquandle $BQ^2_1=\{b, w\}$. The operations are $b^x=w, w^x=b$ and  $b_x=w, w_x=b$ for any $x$. The 1-cells 1 and 2 are cycles.

For 2-cells (squares), $\partial bb=\partial ww=0$ and  $\partial bw=2(w-b)=- \partial wb$.
So $H_1^{\rm BR}=\Z\oplus\Z_2$.

For 3-cells (cubes), $\partial bbb=\partial bbw=-\partial bwb= \partial wbb=bb-ww$.
So $H_2^{\rm BR}=\Z\oplus\Z$.

Assume that $X$ is now a rack and replace $BR$ by $R$.
Let $C_n^{\rm D}(X)$ be the subset of $C_n^{\rm R}(X)$ generated
by $n$-tuples $(x_1, \dots, x_n)$
with $x_{i}=x_{i+1}$ for some $i \in \{1, \dots,n-1\}$ if $n \geq 2$;
otherwise let $C_n^{\rm D}(X)=0$. If $X$ is a quandle, then
$\partial_n(C_n^{\rm D}(X)) \subset C_{n-1}^{\rm D}(X)$ and
$C_\ast^{\rm D}(X) = \{ C_n^{\rm D}(X), \partial_n \}$ is a sub-complex of
$C_\ast^{\rm R}(X)$. 
Put $C_n^{\rm Q}(X) = C_n^{\rm R}(X)/ C_n^{\rm D}(X)$ and 
$C_\ast^{\rm Q}(X) = \{ C_n^{\rm Q}(X), \partial*_n \}$,
where $\partial*_n$ is the induced homomorphism.
We shall follow standard practise and denote all boundary maps by $\partial_n$.

For an abelian group $G$, define the chain and cochain complexes
$$
C_\ast^{\rm W}(X;G) = C_\ast^{\rm W}(X) \otimes G, \quad \partial =
\partial \otimes {\rm id}; \ C^\ast_{\rm W}(X;G) = {\rm Hom}(C_\ast^{\rm
W}(X), G), \quad
 \delta= {\rm Hom}(\partial, {\rm id})
$$
where ${\rm W}$ 
 $={\rm D}$, ${\rm R}$ or ${\rm Q}$.

The $n$\/th {\bf rack/degenerate/quandle homology groups} and the $n$\/th {\bf rack/degenerate /quandle 
cohomology groups\/} can now be defined in the usual way where $R=$ rack, $D=$ degenerate and $Q=$ quandle, see \cite{CJKS}, \cite{NP}, \cite{FRSe}.

For example if $X=Q^3_3$, the 3-colouring quandle then $H_2^Q=0$ and   $H_3^Q=\Z_3$ generated by the cycle $abc+aca$, \cite{NP} \cite{RS}.


\section{How coloured knots determines homology and homotopy classes}
Consider a cellular diagram coloured by a biquandle $X$. Suppose initially that it is edge coloured. Then this defines a map of the underlying surface of the diagram into the classifying space of the biquandle as follows. The diagram defines a cellulation of the surface but we consider the dual cell complex. Around each crossing is a square which is mapped to the named 2-cube in the classifying space. The edges of the square define a normal bundle of the diagram arcs. Each fibre being labelled by the label of the arc. This defines the map on the fibres to the corresponding 1-cube. The points of the regions outside are mapped to the base point. Reidemeister moves correspond to a homotopy of the map. For more details see \cite{FRSa-b}.

For the 3-coloured trefoil the map defines a non-trivial element of $\pi_2$ showing that the trefoil is indeed non-trivial.\footnote{${}^1$}{This is somewhat self referential as having a non-trivial 3-colouring implies a non-zero element of $\pi_2$}

For a whole colouring the construction is similar but defines a 3-cycle. Looking at the trefoils coloured by pairs, since $H_3=\Z_3$ for the 3-colouring space and since the resulting homology class is $+1$ for the right handed trefoil and $-1$ for the left handed trefoil we see that they are distinct.

Whether a diagram can be alternately oriented depends on its class in $H_2$ of the black and white biquandle. 2-sidedness or chessboard colouring depends on its element in $H_3$.

\section{Bibliography}

{\bf B}
A V Bolsinov 
{\it Fomenko invariants in the theory of integrable Hamiltonian systems}
1997 Russ. Math. Surv. 52 997

{\bf BF}
A. Bartholomew and R. Fenn: 
{\it Quaternionic invariants of virtual knots and links}, 
Knot Theory and its Ramifications 17 (2008) 231-251.

{\bf BaF}
A. Bartholomew and R. Fenn: 
{\it Biquandles of small size and some invariants of virtual and welded knots}, 
Knot Theory and its Ramifications 20 (2011) 1-12.

{\bf BuF}
S. Budden and R. Fenn: 
{\it The equation $[B, (A-1)(A,B)]=0$ and virtual knots and links}, 
Fund. Math. {\bf 184} (2004), 19 --29. 

{\bf BuF2}
S. Budden and R. Fenn: 
{\it Quaternion algebras and invariants of virtual knots and links II}, 
JKTR {\bf 17}  (2008) pp ??

{\bf CJKS} 
J. Scott Carter, Daniel Jelsovsky, Seiichi Kamada and Masahico Saito, 
{\it Quandle homology groups, their Betti numbers, and virtual knots}, 
J. Pure Appl. Algebra 157 (2001) 135--155. 

{\bf FJK}
R. Fenn, M. Jordan, and L. Kauffman: 
{\it Biquandles and virtual links}, 
Topology Appl. {\bf 145} (2004), 157--175.

{\bf FRSa}
R. Fenn, C. Rourke and B. Sanderson: 
{\it An introduction to species and the rack space}, 
``Topics in Knot Theory'' (M. E. Bozhuyu, ed.), Kluwer Academic, pp. 33--55, (1993). 

 {\bf FRSb}
R. Fenn, C. Rourke and B. Sanderson: 
{\it Trunks and classifying spaces,}
Appl. Categ. Structures {\bf 3} (1995), 321--356.

{\bf FRSc} 
R. Fenn, C. Rourke and B. Sanderson: 
{\it James bundles and applications,} preprint (1996), 
available at
 http://www.maths.warwick.ac.uk/~cpr/ftp/james.ps

{\bf FRSd} 
R. Fenn, C. Rourke and B. Sanderson: 
{\it The rack space,} TAMS, 359 (2007) 701-740\nl
preprint (2003), 
available
at
arXiv: math.GT/0304228

{\bf FRSe}
R. Fenn, C. Rourke and B. Sanderson: 
{\it James bundles}, 
Proc. London Math. Soc. {\bf 89} (2004), 217--240.  

{\bf FT}
R. Fenn and V. Turaev: 
{\it Weyl Algebras and Knots}, 
J. Geometry and Physics {\bf 57} (2007), 1313-1324

{\bf F}
R. Fenn: 
{\it Quaternion algebras and invariants of virtual knots and links I}, 
JKTR {\bf 17}  (2008) pp 279-304

{\bf History}
http://www.maths.susx.ac.uk/Staff/RAF/Maths/history1.jpg

{\bf IKK}
Atsushi Ishii, Naoko Kamada, Seiichi Kamada:
{\it The virtual magnetic Kauffman bracket skein module and skein relations for the F-polynomial.}
Knot Theory and its Ramifications  17, Issue: 6(2008) pp. 675-688   

{\bf K1} L. H. Kauffman,{\it Virtual Knot Theory }, European J. Comb. (1999) Vol. 20,
663-690.

{\bf K2} L. H. Kauffman, {\it Introduction to virtual knot theory }(arXiv:1101.0665),
in {\it Introductory Lectures on Knot Theory - Selected
Lectures from  the Summer School and Conference on Knot
Theory (ICTP, May 2009)}, edited by L H. Kauffman, S. Lambropoulou, S.
Jablan and J. Przytycki, (to appear in the World Scientific Series on
Knots and Everything (2011)).

{\bf G} Mikhail Goussarov, Michael Polyak and Oleg Viro, {\it Finite type invariants of
classical and virtual knots}, math.GT/9810073.

{\bf Kuper}
G. Kuperberg: 
{\it What is a virtual link?}, 
Algebr. Geom. Topol. {\bf 3} (2003), 587--591. 

{\bf M}
Vassily Olegovich Manturov 
{\it On Virtual Crossing Numbers for Virtual Knots}
arXiv:1107.4828

{\bf YM}
Y. Miyazawa:
{\it Magnetic graphs and an invariant for virtual knots},
Proceedings of Intelligence of Law dimensional Topology
2004, 67--74

{\bf NP}
Maciej Niebrzydowski and Jozef H. Przytycki
{\it Homology of Dihedral Quandles}
Journal of Pure and Applied Algebra
Volume 213, Issue 5, May 2009, Pages 742-755 

{\bf RS}
 C. Rourke and B. Sanderson: 
{\it There are two 2-twist-spun trefoils}
http://msp.warwick.ac.uk/~cpr/

{\bf Saw}
 J. Sawollek:   
{\it On Alexander-Conway polynomials for virtual knots and links\/}, 
J. Knot Theory Ramifications 10 (2001), 151--160
available at
 arXiv: math.GT/9912173

\cite{Stan} { D. Stanovsky} {\it On axioms of biquandles,} J. Knot Theory Ramifications 15/7 (2006) 931--933.

{\bf S}
{ John Stillwell} {\it Classical Topology and Combinatorial Group Theory.} Springer 1993.
\bye


{\bf As} H. Aslaksen: 
{\it Quaternionic determinants}, 
Math. Intel. {\bf 18} (1996), 1-19

{\bf BF}
A. Bartholomew and R. Fenn: 
{\it Quaternionic invariants of virtual knots and links}, 
Knot Theory and its Ramifications 17 (2008) 231-251.

{\bf BaF}
A. Bartholomew and R. Fenn: 
{\it Biquandles of small size and some invariants of virtual and welded knots}, 
Knot Theory and its Ramifications 20 (2011) 1-12.

A. Bartholomew, R. Fenn, N. Kamada, S. Kamada.
New invariants of long virtual knots
Kobe J. Math., 27 (2010) 21–33

{\bf BuF}
S. Budden and R. Fenn: 
{\it The equation $[B, (A-1)(A,B)]=0$ and virtual knots and links}, 
Fund. Math. {\bf 184} (2004), 19 --29. 


{\bf CJKS2001a} 
J. S. Carter, D. Jelsovsky, S. Kamada and M. Saito: 
  {\it Quandle homology groups, their Betti numbers, and virtual knots}, 
J. Pure Appl. Algebra {\bf 157} (2001), 135--155. 


{\bf CKS2001} 
J. S. Carter, S. Kamada and M.  Saito:  
{\it Geometric interpretations of quandle homology},  
J. Knot Theory Ramifications {\bf 10} (2001), 345--386. 

{\bf CKS2002} 
J. S. Carter, S. Kamada and M.  Saito:  
{\it  Stable equivalence of knots on surfaces and virtual knot cobordisms},  
J. Knot Theory Ramifications {\bf 11} (2002), 311--322. 

{\bf CF}
R. H. Crowell and R. H. Fox: 
{\it An introduction to knot theory}, 
Ginn and Co, (1963).

{\bf F}
R. Fenn: 
{\it Quaternion algebras and invariants of virtual knots and links I}, 
JKTR {\bf 17}  (2008) pp 279-304

{\bf FJK}
R. Fenn, M. Jordan, and L. Kauffman: 
{\it Biquandles and virtual links}, 
Topology Appl. {\bf 145} (2004), 157--175.

{\bf FT}
R. Fenn and V. Turaev: 
{\it Weyl Algebras and Knots}, 
J. Geometry and Physics {\bf 57} (2007), 1313-1324

{\bf FR} 
R. Fenn and C. Rourke: 
{\it Racks and links in codimension two}, 
J. Knot Theory Ramifications {\bf 1} (1992), 343--406. 

{\bf FRSa}
R. Fenn, C. Rourke and B. Sanderson: 
{\it An introduction to species and the rack space}, 
``Topics in Knot Theory'' (M. E. Bozhuyu, ed.), Kluwer Academic, pp. 33--55, (1993). 

 {\bf FRSb}
R. Fenn, C. Rourke and B. Sanderson: 
{\it Trunks and classifying spaces,}
Appl. Categ. Structures {\bf 3} (1995), 321--356.

{\bf FRSc} 
R. Fenn, C. Rourke and B. Sanderson: 
{\it James bundles and applications,} preprint (1996), 
available
at
 http://www.maths.warwick.ac.uk/~cpr/ftp/james.ps

{\bf FRSd} 
R. Fenn, C. Rourke and B. Sanderson: 
{\it The rack space,} TAMS, 359 (2007) 701-740\nl
preprint (2003), 
available
at
arXiv: math.GT/0304228

{\bf FRSe}
R. Fenn, C. Rourke and B. Sanderson: 
{\it James bundles}, 
Proc. London Math. Soc. {\bf 89} (2004), 217--240.

{\bf GPV}
M. Goussarov, M. Polyak, and O. Viro: 
{\it Finite-type invariants of classical and virtual knots\/}, 
Topology {\bf 39} (2000), 1045--1068.

{\bf Greene} 
M. T. Greene: 
{\it Some results in geometric topology and geometry,}
Ph.D. Dissertation, Warwick (1997).


{\bf JKS}
F. Jaeger, L. H. Kauffman, and H. Saleur: 
{\it The conway polynomial in $R^3$ and in thickened surfaces: 
A new determinant formulation}, 
J. Combin. Theory Ser. B 
{\bf 61} (1994),  237--259.  





{\bf kamD} N. Kamada:  
{\it The crossing number of alternating link diagrams of a surface\/}, 
Proceedings of Knots 96, World Scientific Publishing Co.,1997, 377--382.


{\bf kamDD}
N.~Kamada:
{\it Span of the Jones polynomial of an alternating virtual link \/},
Algebr. Geom. Topol. {\bf 4} (2004), 1083--1101.


{\bf kamE} 
N.~Kamada: 
{\it A relation of Kauffman's $f$-polynomials of 
virtual links\/}, 
Topology and its Application {\bf 146--147} (2005), 123--132.

{\bf kk} 
N.~Kamada and S.~Kamada: 
{\it Abstract link diagrams and virtual knots\/}, 
J. Knot Theory Ramifications {\bf 9} (2000),
93--106.


 
{\bf SkamA}
S. Kamada:
{\it Braid presentation of virtual knots and welded knots},
Osaka J. Math., to appear,
available
at
 arXiv: math.GT/0008092

{\bf SkamB}
S. Kamada:
{\it Invariants of virtual braids and a remark on left stabilizations
and virtual exchange moves}, 
Kobe J. Math. {\bf 21} (2004), 33--49.





{\bf kauD}
L.~H.~Kauffman: 
{\it Virtual knot theory\/},
Europ.~J.~Combinatorics 
{\bf 20} (1999) 663--690. 


{\bf kaulam}L.~H.~Kauffman and S. Lambropoulou: 
{\it Virtual braids\/}, Fund. Math. {\bf 184} (2004), 159--186. 




{\bf Kuper}
G. Kuperberg: 
{\it What is a virtual link?}, 
Algebr. Geom. Topol. {\bf 3} (2003), 587--591. 


{\bf NP}
Maciej Niebrzydowski and Jozef H. Przytycki
{\it Homology of Dihedral Quandles}
Journal of Pure and Applied Algebra
Volume 213, Issue 5, May 2009, Pages 742-755 

{\bf Saw}
 J. Sawollek:   
{\it On Alexander-Conway polynomials for virtual knots and links\/}, 
J. Knot Theory Ramifications 10 (2001), 151--160
available at
 arXiv: math.GT/9912173

{\bf SWb} D. S. Silver and S. G. Williams,  
{\it Alexander groups and virtual links\/}, 
preprint. 

{\bf SWA} 
D. Silver and S. Williams:
{\it Polynomial invariants of virtual links},
J. Knot Theory Ramifications {\bf 12} (2003), 987--1000. 



\bye
\bye
$$\matrix{%
(pp)&(mm)&(om)&(op)&(po)&(mo)\cr
1&-1&1&-1&0&-1\cr
1&-1&1&-1&1&0\cr
-1&1&-1&1&0&1\cr
-1&1&-1&1&-1&0\cr
1&0&1&0&-1&0\cr
0&1&0&1&0&-1\cr
0&0&0&0&2&-1\cr
0&0&0&0&-1&2\cr
}$$

In this section we show how determinants of matrices with entries in a
general ring $R$ can be defined. Here $R$ is a (possibly
non-commutative), associative ring.  The definition will depend upon a
representation of $R$ into the ring,  $M_{d,d}\Lambda$ of $d\times d$ square matrices with entries in a commutative ring $\Lambda$. For positive codimension determinants it will be useful  if $\Lambda$ supports a greatest common divisor function, written gcd. A good reference for this is \cite{As}.

As an illustration, consider the
non-commutative ring $R = {\H}[t, t^{-1}]$ of Laurent polynomials
whose coefficients are quaternions and the variable $t$ is central. The commutative ring is $\Lambda={\C}[t, t^{-1}]$, the Laurent polynomials in $t$ with complex coefficients. The representation
$\mu: R \to {M}_{2,2}( {\C}[t, t^{-1}])$ has $d=2$ and is
defined by
$$(\alpha_1  + \alpha_2 i + \alpha_3 j + \alpha_4 k) t^s \mapsto 
\pmatrix{(\alpha_1 + \alpha_2 i)  t^s &  (\alpha_3 + \alpha_4 i)   t^s  \cr
(-\alpha_3 + \alpha_4 i)  t^s & (\alpha_1 -\alpha_2 i)  t^s  \cr}
$$
where $\alpha_1, \dots, \alpha_4 \in {\R}$ and $s  \in {\Z}$. 
This representation is {\bf standard} and is usually used in illustrations. 

The ring ${\C}[t, t^{-1}]$ has greatest common divisors, see
\cite{CF}.

Returning now to the general case, let $\eta: R\to M_{d,d}\Lambda$ be
the representation and let $P\in M_{n,m}(R)$. A square submatrix $B$
of $P$ is said to have {\rm codimension} $r$ if it is obtained by
deleting $n-m+r$ rows and $r$ columns if $n\ge m$ or by deleting
$m-n+r$ columns and $r$ rows if $m\ge n$. For simplicity assume $m\ge n$.  Let $B_1, \dots, B_s$ be the codimension $r$ submatrices of $P$ of size $(n-r)\times (n-r)$. Consider $\eta(B_1), \dots, \eta(B_s)$,
which are $d(n-r)\times d(n-r)$ matrices whose entries belong to
$\Lambda$.  The {\it codimension $r$ $\eta$-determinant}, of $P$ is
the gcd of the usual determinants of these matrices.  We denote it by
$\det_\eta^{(r)}(P)$. If $P$ is square then $\det_\eta^{(0)}(P)$ is defined even if $\Lambda$ does not possess greatest common divisors. All determinants are well defined up to multiplication by a
unit.

Now suppose that ${\cal M}$ is a finitely presented $R$-module with
presentation matrix $P\in M_{m,n}(R)$. Then the elements
$\det_\eta^{(r)}(P)$ will be invariants of the module.

\section{The Invariant Modules}

Suppose that $R$ is an associative ring
and
$$S=\pmatrix{A&B\cr C&D\cr}$$
is a $2 \times 2$ matrix with entries from $R$. If $S$ is invertible
and satisfies the set theoretic Yang-Baxter equation in the sense of
\cite{BaF, FJK} then $S$ is called a {\em switch}. The
universal case occurs when $R={\cal F}$ and ${\cal F}$ has the
presentation
$${\cal F}=<A,B\mid A^{-1}B^{-1}AB-B^{-1}AB=BA^{-1}B^{-1}A-A>$$
and $C, D$ are defined by
$$C=A^{-1}B^{-1}A(1-A),\quad D=1-A^{-1}B^{-1}AB.$$

Let $x_0, x_1, \dots, x_n$ be the semi-arcs of a diagram of a long
virtual knot $K$, which appear in this order along $K$.  For each
positive crossing, we consider a relation
$$
S 
\pmatrix{
x_i \cr
x_j\cr}
=
\pmatrix{x_{j+1} \cr
x_{i+1}\cr}
$$
where $x_i$ and $x_j$ are incoming semi-arcs and $x_{j+1}$ and
$x_{i+1}$ are outgoing semi-arcs such that $x_j$ and $x_{j+1}$ are
under-arcs.  For each negative crossing, we consider a relation that
is the inverse of the positive one . 
$$
S 
\pmatrix{
x_{j+1} \cr
x_{i+1}\cr}
=
\pmatrix{x_{i} \cr
x_{j}\cr}
$$
Of course for a virtual crossing the labeling carries over. See the following diagram, (cf. \cite{BaF, FJK}).
\diagram
The corresponding relations are
$$x_{j+1}=Ax_i+Bx_j,\quad x_{i+1}=Cx_i+Dx_j$$
and
$$x_{i}=Ax_{j+1}+Bx_{i+1},\quad x_{j}=Cx_{j+1}+Dx_{i+1}$$
The module ${\cal M}_K$ is the $R$-module generated by $x_0, x_1,
\dots, x_n$ with the relations associated with positive crossings and
negative crossings.  There is one more generator than relation.

The module $\widehat{\cal M}_K$ is the quotient of ${\cal M}_K$ by an
additional relation $x_0=x_n$. The number of generators and relations
are the same. That is, the presentation matrix is square.

The module ${\oM}_K$ is the quotient of ${\cal M}_K$ by an
additional relation $x_0=0$. Again, the presentation matrix is square.

The module ${\nM}_K$ is the quotient of ${\cal M}_K$ by an
additional relation $x_n=0$. The generator $x_n$ is the label on the outgoing arc. Again, the presentation matrix is square.

\theorem{
Suppose $A, 1-A$ and $B$ are invertible. Then the modules ${\cal M}_K$,
$\widehat{\cal M}_K$, ${\oM}_K$ and  ${\nM}_K$ are invariants of a long virtual knot $K$.}

{\it Proof}: Since $S$ is invertible, satisfies the set theoretic
Yang-Baxter equation, and since $1-A$ is invertible these module are
preserved by all generalized Reidemeister moves (see \cite{BaF,
FJK}).  \qed

As an illustration consider the ``fly'' long knot, $F$,  pictured below.

\diagram

The presentation matrices in the four cases are
$${\cal M}=\pmatrix{
-1&B&0&A&0\cr
0&D&-1&C&0\cr
0&A&0&B&-1\cr
0&C&-1&D&0\cr},\
\widehat{\cal M}=\pmatrix{
-1&B&0&A&0\cr
0&D&-1&C&0\cr
0&A&0&B&-1\cr
0&C&-1&D&0\cr
1&0&0&0&-1\cr}$$
and
$$
{\oM}=\pmatrix{
-1&B&0&A&0\cr
0&D&-1&C&0\cr
0&A&0&B&-1\cr
0&C&-1&D&0\cr
1&0&0&0&0\cr},\
{\nM}=\pmatrix{
-1&B&0&A&0\cr
0&D&-1&C&0\cr
0&A&0&B&-1\cr
0&C&-1&D&0\cr
0&0&0&0&1\cr},\
$$
\section{The Invariant Polynomials}

Let $K$ be a long knot and let ${\cal M}_K$, $\widehat{\cal M}_K$,
${\oM}_K$  and ${\nM}_K$ be the ${\cal F}$-modules defined in the previous section. Let $P$,
$\widehat P$,  $\oP$ and $\nP$ be the respective presentation
matrices. Suppose we now represent the algebra as matrices so that we can define determinantal invariants as described in section 2. Each entry in the $P$
matrices is a $d\times d$ matrix for some $d$. Let $p^{(r)}_K$,
$\widehat{p}^{(r)}_K$, 
$\op^{\raise -8pt\hbox{$\scriptstyle(r)$}}_K$ 
and 
$\np^{\raise -8pt\hbox{$\scriptstyle(r)$}}_K$ 
be the corresponding
determinants in the commutative ring $\Lambda$, with codimension $r=0,1,2\ldots$.

Let $\widehat{K}$ be the closure of the long knot $K$. Then
$\widehat{K}$ also has an invariant ${\cal F}$-module, see \cite{F}.
Let $q^{(r)}_{\widehat{K}}$ be the sequence of determinants in $\Lambda$,
corresponding to the presentation.

\theorem{
\parindent=20pt
\item{ (1)} $q^{(r)}_{\widehat{K}}
=\widehat{p}^{(r)}_K$ 
\item{ (2)} $p^{(r)}_K$ divides $\op^{\raise -8pt\hbox{$\scriptstyle(r)$}}_K$. 
\item{ (3)} $p^{(r)}_K$ divides $\np^{\raise -8pt\hbox{$\scriptstyle(r)$}}_K$. 
}
\parindent=0pt

{\it Proof}:
\def\oQ{\buildrel o \over {Q}}

 
(1) This follows since $\widehat{\cal M}_K$ is
equal to the module of
the closure $\widehat{K}$ of $K$.
 
(2) The module ${\cal M}_K$ has an $n \times (n+1)$
matrix $P$ as a
presentation matrix such that the first column corresponds to $x_0$
and the last column corresponds to $x_n$. A codimension $r$ submatrix
is obtained from $P$ by deleting $r$ rows and $r+1$ columns. Let $B_1,
\dots, B_s$ be the codimension $r$ submatrices of $P$. Then
$p^{(r)}_K$ divides the determinant of all of these after the
representation and is the largest, by division, element which does so.
 
The module ${\oM}_K$ has an $(n+1) \times (n+1)$
presentation matrix $\oP$ that is obtained from $P$ by adding
the row $(1, 0, \dots, 0)$ to the bottom. 
The presentation matrix $\oP$ is simplified to $\oQ$ that is an
$n \times n$ matrix obtained from $\oP$ by deleting the first column  and the bottom row, which  is obtained from $P$ by deleting the first column.
Since $\oQ$ is
square a codimension $r$ submatrix is obtained from $\oQ$ by
deleting $r$ rows and $r$ columns. It follows that
$\op^{\raise -8pt\hbox{$\scriptstyle(r)$}}_K$
 is the gcd of a subset of the values for which
$p^{(r)}_K$ is the gcd. The result now follows.
 
(3) The proof is similar to (2). \qed
 

\theorem{ Let $K_1 \cdot K_2$ be the concatenation product of
two long virtual knots $K_1$ and $K_2$. For any suitable
representation of the fundamental modules
$$\op^{\raise -8pt\hbox{$\scriptstyle(0)$}}(K_1 \cdot K_2) = 
\op^{\raise -8pt\hbox{$\scriptstyle(0)$}}(K_1) 
\op^{\raise -8pt\hbox{$\scriptstyle(0)$}}(K_2)$$
and
$$\np^{\raise -8pt\hbox{$\scriptstyle(0)$}}(K_1 \cdot K_2) = 
\np^{\raise -8pt\hbox{$\scriptstyle(0)$}}(K_1) 
\np^{\raise -8pt\hbox{$\scriptstyle(0)$}}(K_2). $$
}

{\it Proof}: Let $P_1= (a_0 a_1 \cdots a_n)$ and $P_2 =(b_0 b_1 \cdots
b_{n'})$ be the presentation matrices of ${\cal M}_{K_1}$ and ${\cal
M}_{K_2}$ associated with their diagrams.  Then ${\oM}_{K_1}$,
${\oM}_{K_2}$ and ${\oM}_{K_1 \cdot K_2}$ have presentation matrices
$$
 (a_1  \cdots  a_n), \quad  (b_1 \cdots  b_{n'}), \quad 
\pmatrix{
a_1  \cdots  a_{n-1}  &      a_n        & 0 \cdots 0 \cr
  0 \cdots   0              &           b_0  &   b_1  \cdots   b_{n'} \cr}
$$
respectively.  Thus we have the result for $\op$. The proof for $\np$ is
similar.
\qed
\section{Simplifying the Modules and some Calculations}
The presentation matrices defined above can be simplified by the usual rules for manipulating non-commuting relations. That is
{\parindent=20pt
\item{1} Interchange any row(column).
\item{2} Multiply any row(column) on the left(right) by a unit.
\item{3} Add any row(column) multiplied on the left(right) to a different row(column).
\item{4} Introduce or delete any zero row.
\item{5} $P\leftrightarrow \pmatrix{1&0\cr 0& P\cr}$
}

If we apply these rules to the example in section 3 and simplify as far as possible we get,
$$\pmatrix{
D-C&D-C\cr},\
\pmatrix{0\cr},\
\pmatrix{(C-D)(1+B^{-1}A)\cr},\
\pmatrix{(D-C)(1+B^{-1}A)\cr}$$
for the four presentation matrices. Note that the presentation matrix for $\widehat{\cal M}_K$ will reduce to zero since the closure of $B$ is the trivial knot.

Now apply the homomorphism which replaces $A, B, C, D$ with $1+i, -tj, t^{-1} j, 1+i$ respectively and use the standard representation. The three codimension zero polynomial invariants are
$|D-C|^2=2+t^{-2}$, $0$ and $|(C-D)(1+B^{-1}A)|^2=(2+t^{-2})(1+2t^{-2})$. All the higher codimension polynomials are 1.
\section{Symmetries of Long Virtual Knots}
There are various symmetries of the knot diagram which can be applied.
Consider reflection in the plane of a knot diagram $D$. Let $-D$ denote the resulting diagram. This interchanges plus and minus crossings. Let $\overline D$ denote the result of reflection in the $x$-axis. Finally let $D^*$ be obtained by reversing the arrow and rotating the result through 180 degrees. 

For the fly $F^*=\overline F$.
The effect of the other three possibilities on the fly are illustrated below.
\medskip
\centerline{\ddiagram\quad \ddiagram\quad \ddiagram}
\centerline{\npt The fly reflected: $-F$,\quad $\overline F$,\quad $-\overline F$}
Using a suitable representation of the fundamental algebra all three can be distinguished from themselves and the original fly as follows.

The resulting polynomials are tabulated as follows. The switch used is given by a representation of the quantum Weyl
algebra, with  $A, B, C, D$ the following $2\times2$ matrices.

$$
	A = \pmatrix{1-q & -q^3+2q^2-1 \cr 0 & 1-q },\ 
	B = \pmatrix{q & 1 \cr 0 & 1} 
$$
$$
	C = \pmatrix{1 & (-q^4+3q^3-2q^2-2q+1)/q \cr 0 & q},\ 
	D = \pmatrix{0 & (q^3-2q^2+1)/q \cr 0 & 0 } 
$$

$$\vbox{ \offinterlineskip \halign{
     \vrule # & \hfil \quad $#$ \quad \hfil & 
     \vrule # & \hfil     \ $#$ \     \hfil & 
     \vrule # & \hfil     \ $#$ \     \hfil & 
     \vrule # & \hfil     \ $#$ \     \hfil &
     \vrule # & \hfil     \ $#$ \     \hfil &
	 \vrule#\cr 
\noalign{\hrule} 
& K && p^{(0)}(K) && \widehat{p}^{(0)}(K) && \op^{\raise
-8pt\hbox{$\scriptstyle(0)$}}(K) && \np^{\raise
-8pt\hbox{$\scriptstyle(0)$}}(K) & \cr 
\noalign {\hrule} 
height1pt&\omit&height1pt&\omit&height1pt&\omit&height1pt&\omit&height1pt&\omit&height1pt\cr
& F           && 1 && 0 && (2-q)/q && (2-q)/q & \cr 
\noalign {\hrule} 
height1pt&\omit&height1pt&\omit&height1pt&\omit&height1pt&\omit&height1pt&\omit&height1pt\cr
&-F           && 2-q && 0 && (2-q)/q && (2-q)/q & \cr 
\noalign {\hrule} 
height1pt&\omit&height1pt&\omit&height1pt&\omit&height1pt&\omit&height1pt&\omit&height1pt\cr
&\overline F  && 1  && 0 && 2q-1 && 2q-1 & \cr
\noalign {\hrule} 
height1pt&\omit&height1pt&\omit&height1pt&\omit&height1pt&\omit&height1pt&\omit&height1pt\cr
&-\overline F &&  2q-1 && 0 && 2q-1 && 2q-1 &\cr 
\noalign {\hrule} 
}}$$


\theorem{Consider the following three conditions on a switch $S$.
a) $S=S^{\dag}$, b) $S^2=1$, c) $SS^{\dag}=1$, where
$S^{\dag}=\pmatrix{D&C\cr B&A\cr}$.
\nl
If a) then $A=D$ and $B=C$ and $K$ cannot be distinguished from $-\overline K$. \nl
If b) then the underlying algebra is the Weyl algebra and $K$ cannot be distinguished from $-K$ or $K^*$. \nl
If c) then $A, B$ commute and $S$ is $\pmatrix{2&\pm1\cr \mp1&0\cr}$, a specialization of the Alexander switch. Moreover $K$ cannot be distinguished from $\overline K$.}

{\it Proof}: Most of the results easily follow by looking at the conditions on the entries of $S$ and how this affects the calculations of the modules.
\nl
For b) the underlying algebra is the Weyl algebra because of results in \cite{FT}.
\nl
For c) the condition $SS^{\dag}=1$ implies 
$$\matrix{
AD+B^2&=\hfill1,&AC&=-BA \cr CD&=-DB,&C^2+DA&=\hfill1\cr}.$$
Using the relations
$$C=A^{-1}B^{-1}A(1-A),\quad D=1-A^{-1}B^{-1}AB$$
gives the result. \qed 

It is well known that the product of two classical knots is a commutative operation. This is not the case for the product of two long virtual knots.
For example the products $F\cdot \overline F$ and $\overline F\cdot F$ are distinct. 
A calculation using the Budapest switch shows that for $F\cdot \overline F$, $ 
p^{(0)} = 6t^4+15t^2+6$ whereas for $\overline F\cdot F$ we have
$p^{(0)} = 3t^4+15/2t^2+3$, which is half the previous polynomial. If we are working over the integer quaternions then 2 is not a unit and so this shows the knots are distinct. 

This is not perhaps a ``killer" example. If we consider
$$(F \cdot \overline{F}) \cdot F \hbox{ and  }F\cdot (F \cdot \overline{F})$$
then in both cases $p^{(0)} =-12t^8-60t^6-99t^4-60t^2-12$ but the values of $p^{(1)} $ are $9(t^2+1)$ and $-9/2(2t^2+1)(t^2+2)$ respectively.

\section{Long flat virtual knots}  

We now repeat the previous analysis for long flat knots. For a full
discusion of the following see \cite{FT}. A switch $S$ can be used
provided it satisfies $S^2=id$. The conditions for this are contained
in the following theorem.

\theorem{
Suppose $A, 1-A$ and $B$ are invertible and
$$
S = 
\pmatrix{
A& B \cr
C& D\cr
}$$
is a $2 \times 2$ matrix with entries satisfying
$$A^{-1}B^{-1}AB-B^{-1}AB=BA^{-1}B^{-1}A-A=1$$
and $C, D$ are defined by
$$C=A^{-1}B^{-1}A(1-A),\quad D=1-A^{-1}B^{-1}AB.$$
Then $S^2=1$ if and only if $u=B, v=B^{-1}A^{-1}$ satisfy $uv-vu=1$.
(The elements $u, v$ are the generators of the Weyl algebra.)}
\qed

We now repeat the construction considered earlier and arrive at modules ${\cal WM}_F$, $\widehat{\cal WM}_F$, ${\cal W}\!\oM_F$ and ${\cal W}\!\nM_F$ which are invariants of a long flat knot $F$.

Again by analogy we can find tractible invariants given representations
onto finite matrices. Many examples are given in \cite{FT}.

The following is an example with ring ${\Z}_2[a, x, y]$.
$$
u = \pmatrix{
x  &  a \cr
0  &  x \cr},\
v = \pmatrix{
y  &  0 \cr
1/a  &  y \cr}
$$
Consider the ``flat fly" denoted by $FF$ and its reflection in the $x$-axis illustrated below.
\medskip
\centerline{\ddiagram,\quad\ddiagram}

\centerline{$FF$ and $\overline{FF}$}

Using the representation above, setting $y=x$ and $a=1$, the codimension zero polynomials are

$$\vbox{ \offinterlineskip \halign{
     \vrule # & \hfil \quad $#$ \quad \hfil & 
     \vrule # & \hfil     \ $#$ \     \hfil & 
     \vrule # & \hfil     \ $#$ \     \hfil & 
     \vrule # & \hfil     \ $#$ \     \hfil &
     \vrule # & \hfil     \ $#$ \     \hfil &
	 \vrule#\cr 
\noalign{\hrule} 
& K && p^{(0)}(K) && \widehat{p}^{(0)}(K) && \op^{\raise
-8pt\hbox{$\scriptstyle(0)$}}(K) && \np^{\raise
-8pt\hbox{$\scriptstyle(0)$}}(K) & \cr 
\noalign {\hrule} 
height1pt&\omit&height1pt&\omit&height1pt&\omit&height1pt&\omit&height1pt&\omit&height1pt\cr
& {\bf d}(FF)  && x^2+1 && 0 && x^8+x^6+x^2+1 && x^8+x^6+x^2+1 & \cr 
\noalign {\hrule} 
height1pt&\omit&height1pt&\omit&height1pt&\omit&height1pt&\omit&height1pt&\omit&height1pt\cr
&{\bf d}(\overline{FF}) && x^6+1 && 0 && x^8+x^6+x^2+1 && x^8+x^6+x^2+1 & \cr 
\noalign {\hrule} 
}}$$

This shows that $FF$ and $\overline{FF}$ are both non-trivial and distinct.

It is interesting to note Turaev's descent map ${\bf d}$ of
long flat knots to long virtual knots as an alternative method to show that $FF$ is non-trivial. This lifts flat knots by turning the first time a crossing is met to an overcrossing. For example $FF$ and $\overline{FF}$ are converted as shown in the following diagram.

$$\matrix{\ddiagram\ &\raise 20 pt\hbox{$\to\ $} &\ddiagram\cr}$$
$$\matrix{\ddiagram\ &\raise 20 pt\hbox{$\to\ $} &\ddiagram\cr}$$

Then, using the same switch as in the previous example, the polynomials for ${\bf d}(FF)$ and ${\bf d}(\overline{FF})$ 
are also

$$\vbox{ \offinterlineskip \halign{
     \vrule # & \hfil \quad $#$ \quad \hfil & 
     \vrule # & \hfil     \ $#$ \     \hfil & 
     \vrule # & \hfil     \ $#$ \     \hfil & 
     \vrule # & \hfil     \ $#$ \     \hfil &
     \vrule # & \hfil     \ $#$ \     \hfil &
	 \vrule#\cr 
\noalign{\hrule} 
& K && p^{(0)}(K) && \widehat{p}^{(0)}(K) && \op^{\raise
-8pt\hbox{$\scriptstyle(0)$}}(K) && \np^{\raise
-8pt\hbox{$\scriptstyle(0)$}}(K) & \cr 
\noalign {\hrule} 
height1pt&\omit&height1pt&\omit&height1pt&\omit&height1pt&\omit&height1pt&\omit&height1pt\cr
& {\bf d}(FF)  && x^2+1 && 0 && x^8+x^6+x^2+1 && x^8+x^6+x^2+1 & \cr 
\noalign {\hrule} 
height1pt&\omit&height1pt&\omit&height1pt&\omit&height1pt&\omit&height1pt&\omit&height1pt\cr
&{\bf d}(\overline{FF}) && x^6+1 && 0 && x^8+x^6+x^2+1 && x^8+x^6+x^2+1 & \cr 
\noalign {\hrule} 
}}$$

Our invariants are consequence of biquandles \cite{FJK}.  When we deform a given virtual knot or a given long virtual knot into a braid form, it is easier to calculate the biquandle, the quaternionic module and the determinant invariants.  For braiding of virtual knots and long virtual knots, refer to 
\cite{SkamA, SkamB, kaulam}.  
\section{References}

\bye

\end{thebibliography}
\bye